\documentclass[fleqn]{amsart}
\usepackage{amssymb, graphics, amsmath, epsfig, changebar, enumitem, hyperref}
\usepackage{mathtools} 
\usepackage[percent]{overpic}
\usepackage[caption=false]{subfig} 

\theoremstyle{plain}
\newtheorem{theorem}{Theorem}
\newtheorem{corollary}{Corollary}

\newtheorem{lemma}{Lemma}
\newtheorem{proposition}{Proposition}
\newtheorem{remark}{Remark}

\theoremstyle{definition}

\theoremstyle{example}

\theoremstyle{remark}

\numberwithin{equation}{section}
\newcommand{\N}{\mathbb{N}}
\newcommand{\Z}{\mathbb{Z}}

\newcommand{\Sum}{\displaystyle\sum}

\newcommand{\Om}{\Omega}

\graphicspath{{./images/}}

\newcommand{\overimage}[3]{\begin{overpic}{#1}\put(#2){#3}\end{overpic}}
\newcommand{\oversimage}[5]{\begin{overpic}{#1}\put(#2){#3}\put(#4){#5}\end{overpic}}

\begin{document}                   
\title{On some moves on links and the Hopf crossing number}
\author{Maciej Mroczkowski}
\address{Institute of Mathematics\\
Faculty of Mathematics, Physics and Informatics\\
University of Gdansk, 80-308 Gdansk, Poland\\
e-mail: maciej.mroczkowski@ug.edu.pl}

\begin{abstract}
We consider arrow diagrams of links in $S^3$ and define $k$-moves on such diagrams, for any $k\in\N$. We study the equivalence classes of links in $S^3$
up to $k$-moves. For $k=2$, we show that any two knots are equivalent, whereas it is not true for links.
We show that the Jones polynomial at a $k$-th primitive root of unity is unchanged by a $k$-move, when $k$ is odd. It is multiplied by $-1$, when $k$ is even. 
It follows that, for any $k\ge 5$, there are infinitely many classes of knots modulo $k$-moves.
We use these results to study the Hopf crossing number. In particular, we show that it is unbounded for some families of knots.
We also interpret $k$-moves as some identifications between links in different lens spaces $L_{p,1}$.
\end{abstract}
\maketitle

\let\thefootnote\relax

\section{Introduction}
Several types of moves on links in $S^3$ were studied extensively for a long time (see for example \cite{Fo},\cite{Ptk}).
There are two important questions related to a chosen set of moves $\mathcal M$.
The first is: are any two links $L$ and $L'$ related by moves from $\mathcal M$?
In other words, is it possible to go from $L$ to $L'$ by a finite sequence of such moves and isotopies?
The second question is: given $L$ and $L'$ that are related by moves from $\mathcal M$, what is the least number of such moves needed
to go from $L$ to $L'$?

We introduce $k$-moves, $k\in\N$, on links in $S^3$. They are defined with the help of arrow diagrams of links (see \cite{MD},\cite{M}).
We answer the first of the above questions for $k\neq 3,4$. For $k=1$, all links are related by $k$-moves. For $k=2$ all knots (but not all links)
are related by $k$-moves. For $k\ge 5$, we show that the value of the Jones polynomial in a $k$-th primitive root of unity is unchanged under a $k$-move,
if $k$ is odd, and changes sign, if $k$ is even (see Theorem~\ref{thm:modJones}).
It follows that there are infinitely many classes of knots modulo $k$-moves, when $k\ge 5$.

Let $p:S^3\to S^2$ be the Hopf map and $L\subset S^3$ a link. The minimal number of crossings in $p(L)$ among generic projections
of $L$ in $S^2$ is called the {\it Hopf crossing number} and is denoted $h(L)$. The idea of such and invariant was mentioned by Fiedler in \cite{F}.
In \cite{M}, knots with Hopf crossing number at most $1$ were classified. It was shown there, that $h(L)$ equals the minimum of crossings
among arrow diagrams of $L$. It is also equal to the minimal number of crossings among Turaev's gleams of $L$ (see~\cite{T},\cite{B}).

The introduction of $k$-moves is motivated by the study of the Hopf crossing number. One problem with this invariant is that, for a given integer number $n\ge 0$, 
there a inifinitely many knots $K$, such that $h(K)\le n$. Thus, for example, in contrast to the classical crossing number, 
one cannot list all knots with $h(K)\le n$, in order to show that a knot $K'$,
not in this list, satisfies $h(K')>n$. On the other hand, for a given $k\in\N$, there are finitely many knots $K$ {\it up to $k$-moves}, such that $h(K)\le n$. 
Using this idea and Theorem~\ref{thm:modJones}, we study the Hopf crossing number of some links.
In particular, we exhibit some families of knots with unbounded Hopf crossing numbers, such as some prime alternating knots and
prime knots with braid index~$3$.

In section~\ref{sec:arrow}, we recall the definition of arrows diagrams. In section~\ref{sec:equiv}, we introduce $k$-moves and give an interpretation
of these moves in terms of classical diagrams; then we study the equivalence under $1$ and $2$-moves. In section~\ref{sec:jones}, we prove the main result,
Theorem~\ref{thm:modJones} and some corollaries, establishing the relationship between $k$-moves and values of the Jones polynomial at primitive $k$-roots of unity.
Finally, in section~\ref{sec:appl}, we use the previous results to study $k$-equivalance classes of torus knots $T(n,n+1)$ and exhibit some families of knots
with unbounded Hopf crossing number; we end with an interpretation of $k$-moves as some identifications between links in different lens spaces $L_{p,1}$.

\section{Arrow diagrams of links in $S^3$}\label{sec:arrow}
Arrow diagrams for links in $F\times S^1$, where $F$ is an orientable surface, were introduced in \cite{MD}.
They are like classical diagrams (generically immersed closed curves in $F$ together with information of over/under at the crossings)
with some arrows outside crossings. To obtain an arrow diagram of a link $L$ in $F\times S^1$, one cuts $F\times S^1$ along $F\times\{1\}$
and projects $L$ from $F\times [0,1]$ onto $F$. The arrows denote places where $L$ was cut and point to the part of $L$ closest to $F\times \{0\}$
in $F\times [0,1]$.

\begin{figure}[h]
\centering
\includegraphics{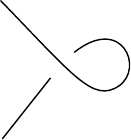} \raisebox{1.9 em}{$\; \overset{\Omega_1}{\longleftrightarrow} \;$} \includegraphics{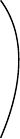}
\hspace{2em}
\includegraphics{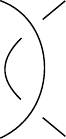} \raisebox{1.9 em}{$\; \overset{\Omega_2}{\longleftrightarrow} \;$} \includegraphics{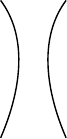}
\hspace{2em}
\includegraphics{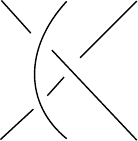} \raisebox{1.9 em}{$\; \overset{\Omega_3}{\longleftrightarrow} \;$} \includegraphics{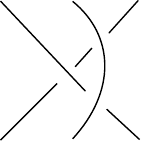}
\hspace{2em} \\[2em]
\includegraphics{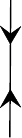} \raisebox{1.7 em}{$\; \overset{\Omega_4}{\longleftrightarrow} \;$} \includegraphics{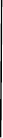}
\raisebox{1.7 em}{$\; \overset{\Omega_4}{\longleftrightarrow} \;$} \includegraphics{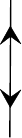} \hspace{2em}
\includegraphics{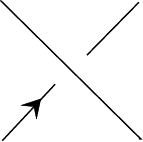} \raisebox{1.7 em}{$\; \overset{\Omega_5}{\longleftrightarrow} \;$} \includegraphics{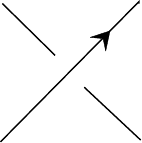}
\\[2em]
\includegraphics{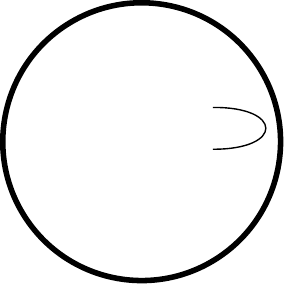} \raisebox{4.2 em}{$ \;
\overset{\Omega_{\infty}}{\longleftrightarrow} \; $ }
\begin{overpic}{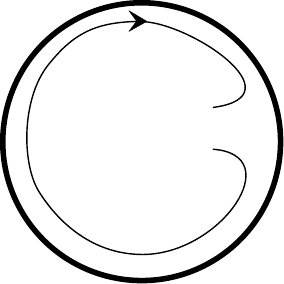}   \end{overpic}
\caption{Reidemeister moves}
\label{reid_moves}
\end{figure}

In \cite{MD}, it was shown that $D$ and $D'$ are two arrow diagrams of the same link in $F\times S^1$, if and only if one can go from $D$ to $D'$ 
with a series of five Reidemeister moves, three classical and two extra ones ($\Om_4$ and $\Om_5$).
When $F$ is a disk $B$, so that $B\times S^1$ is a solid torus $T$, the arrow diagrams lie in $B$.
$S^3$ is obtained by gluing another torus $T'$ to $T$, with the boundary of the meridional disk of $T'$ attached to a $(1,1)$ curve in $\partial T$.
Any link in $S^3=T\cup T'$ can be pushed into $T$, so that is has an arrow diagram in $D$. Gliding an arc through the meridional disk of $T'$ gives
rise to one extra Reidemeister move, denoted $\Om_\infty$. All Reidemeister moves for arrow diagrams of links in $S^3$ are shown in Figure~\ref{reid_moves}.

One can view $T$ and $T'$ as consisting of fibers of the Hopf fibration of $S^3$. Then the projection from $T$ onto $B$ is a restriction of the
Hopf map $p:S^3\to S^2$. A link $L$ in $S^3$ can be isotoped so that it lies in $T$ without changing the number of crossings in $p(L)$. 
Thus, the Hopf crossing number of $L$, or $h(L)$, is the minimum of crossings among all arrow diagrams of $L$. As classical diagrams (without arrows) 
form a subset of arrow diagrams, it is clear that $h(L)\le c(L)$, where $c(L)$ is the crossing number of $L$.
For details see \cite{M}. 

\begin{figure}[h]
\centering
\includegraphics{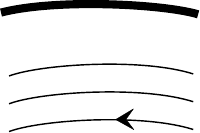} \raisebox{1.9 em}{$\; \overset{\Omega_2,\;\Omega_5}{\xrightarrow{\hspace*{2.5em}}} \;$} \includegraphics{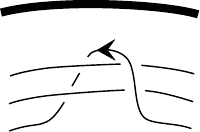}
\raisebox{1.9 em}{$\; \overset{\Omega_1,\;\Omega_5}{\xrightarrow{\hspace*{2.5em}}} \;$} \includegraphics{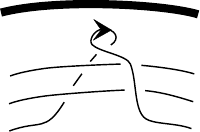}
\raisebox{1.9 em}{$\; \overset{\Omega_\infty,\;\Omega_4}{\xrightarrow{\hspace*{2.5em}}} \;$} \includegraphics{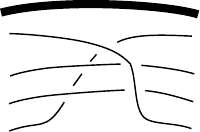}
\caption{Eliminating an arrow}
\label{elim_arr}
\end{figure}

Any arrow can be eliminated with Reidemeister moves, see Figure~\ref{elim_arr} (if the arrow points in the opposite direction, the middle step involving
$\Om_1$ and $\Om_5$ is skipped). 

\section{Equivalence modulo $k$ on links}\label{sec:equiv}
\begin{figure}[h]
\centering
\includegraphics{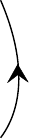}\raisebox{1.5 em}{$\;3$}
\raisebox{1.9 em}{$\;=\;\;\;$}
\includegraphics{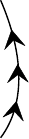}\raisebox{1.9 em}{$\; \overset{3-move}{\xleftrightarrow{\hspace*{2.5em}}} \;$}
\includegraphics{R1b} \raisebox{1.9 em}{$\; \overset{3-move}{\xleftrightarrow{\hspace*{2.5em}}} \;$} \includegraphics{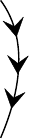}
\raisebox{1.9 em}{$\;=\;\;\;$}\includegraphics{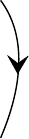}\raisebox{1.9 em}{$\;3$}
\raisebox{1.9 em}{$\;=\;\;\;$}\includegraphics{R1b_k_arr2u}\raisebox{1.5 em}{$\;-3$}
\caption{Examples of $3$-moves}
\label{k_moves}
\end{figure}

Let $D$ be an arrow diagram. For $k\in\N$, a {\it $k$-move} on $D$ is the addition or removal of $k$ arrows on a strand of $D$ without crossings between them,
pointing in the same direction. In Figure~\ref{k_moves}, $3$-moves are shown. By convention, the $3$ next to the arrow stands for $3$ arrows in
the same direction. Also, a negative integer $b$ next to an arrow stands for $|b|$ arrows opposite to the one that is pictured.

We say that two diagrams $D$ and $D'$ are {\it $k$-equivalent}, if there is a series of $k$-moves and Reidemeister moves allowing to pass from $D$ to $D'$.
Two links $L$, with diagram $D$, and $L'$, with diagram $D'$, are {\it $k$-equivalent}, if $D$ and $D'$ are $k$-equivalent.
We also say that $L$ and $L'$ are {\it equivalent modulo $k$}.

\begin{figure}[h]
\centering
\includegraphics{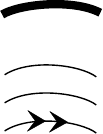}\raisebox{1.9 em}{$\longrightarrow$}\includegraphics{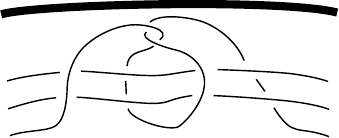}\qquad\qquad
\includegraphics{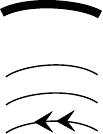}\raisebox{1.9 em}{$\longrightarrow$}\includegraphics{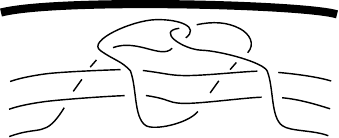}
\caption{Eliminating two consecutive arrows}
\label{2_moves}
\end{figure}

A $k$-move could be defined for classical diagrams (without reference to arrow diagrams). Indeed, if a $k$-move from $D$ to $D'$ adds $k$ arrows,
we could eliminate all other arrows on $D$ and $D'$ (see Figure~\ref{elim_arr}), so that $D$ has no arrows and $D'$ has $k$ arrows. Then, the $k$
arrows on $D'$ could be eliminated as well, yielding a move between classical diagrams. However, in such a setting,
the $k$-move would be more complicated to define. For example, for $k=2$, two cases of eliminating two consecutive arrows, assuming there
are no other arrows, are shown in Figure~\ref{2_moves}. In that figure, there are two strands between the arrows and the boundary of the diagram, 
but in general it could be any number of strands, including zero. Also, we could limit ourselves to just one (anyone) of the two cases shown in that figure.
For example, the two arrows that are pushed to be {\it clockwise} with respect to the boundary (first case in the figure), 
become {\it counterclockwise}, if pushed to the antipodal side of the boundary (second case).

\begin{figure}[h]
\centering
\includegraphics{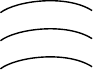}\raisebox{1.9 em}{$\longrightarrow$}\includegraphics{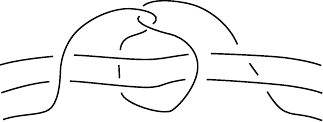}\qquad\qquad
\includegraphics{2mv_c}\raisebox{1.9 em}{$\longrightarrow$}\includegraphics{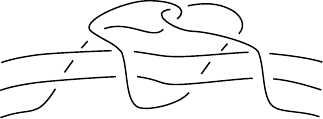}
\caption{$2$-moves for classical diagrams}
\label{2_moves_cl}
\end{figure}

The $2$-moves for classical diagrams are shown in Figure~\ref{2_moves_cl}.
Again, there could be any number of strands in this figure. Notice, that if there are no middle strands, the first case of the move yields an isotopic link,
whereas the second yields a connected sum with a right handed trefoil.

Let $T_n$ be an oval with $n$ counterclockwise arrows on it, if $n\ge 0$, and $|n|$ clockwise arrows, if $n<0$.
It was shown in \cite{M}, that $T_n$ is the right handed torus knot $T(n,n+1)$, for $n\ge0$.
One checks easily that, for $n<0$, $T_n$ can be transformed into $T_{-n-1}$ with  $\Om_\infty$ and $\Om_4$.
For a given $n$, the knots $T_{n+ks}$, $s\in\Z$, are obviously related by $k$-moves. For instance, taking $k=5$, $T_7$ is $5$-equivalent to $T_2$,
$T_{-3}$ (which equals $T_2$), $T_{-8}$ (which equals $T_7$), $T_{12}$, etc. On the other hand, as $T_{-2}$, $T_{-1}$, $T_0$ and $T_1$ are all
trivial, $T_3$, $T_4$, $T_5$ and $T_6$ are $5$-equivalent to the unknot. We will see later, that for the knots $T_n$ there are exactly two equivalence
classes modulo $5$ (see Proposition~\ref{pro:tn}).
 
The knots $T_n$ play a fundamental part for studying $k$-moves. For example, let $K$ be a knot represented by $D$, a diagram without arrows.
Adding $k>0$ arrows on a strand of $D$ next to the boundary of the diagram, we obtain the connected sum $K\sharp T_k$ or $K\sharp T_{-k}=K\sharp T_{k-1}$,
depending on the direction of the $k$ arrows. This follows from the fact that we can eliminate the $k$ arrows keeping freezed all $D$ except for the strand
on which we added these arrows. 
The $k$-move is no longer necessarily a connected sum with $T_k$ or $T_{k-1}$, if we add the arrows on a strand of $D$ which is not next to the boundary, 
or if $D$ has some arrows on it. 

In the rest of this section, we study the question of when are two links or knots equivalent (or not) modulo $k$ (mostly for $k=1$ and $2$).

\begin{remark}
If $D$ is a classical diagram of a link $L$ (i.e. a diagram without arrows) it is possible to change any crossing of $D$ with $1$-moves and $\Om_5$:
with a $1$-move create an arrow pointing to the undercrossing, use $\Om_5$ and remove the arrow with another $1$-move.
Thus any link is $1$-equivalent to a trivial link.
\end{remark}

We show now, that for $k\ge 2$, $k$-moves do not trivialize all links.

\begin{proposition}
The linking number modulo $k$ between any two components of a link is preserved under $k$-equivalence. Thus, for $k\ge 2$, there are
links that are not $k$-equivalent to trivial links.
\begin{proof}
Let $D$ be an arrow diagram and $D'$ be obtained from $D$ by performing a $k$-move. With Reidemeister moves eliminate in the same way
all arrows in $D$ and $D'$ (see Figure~\ref{elim_arr}), except the $k$ arrows created by the $k$-move. Thus, we may assume that $D$ has no arrows
and $D'$ has $k$ arrows all on the same component, say $C_1$. Let $C_2$ be another component.
The $k$ arrows on $D'$ are separated from the boundary of the diagram with some arcs, some of which may belong to $C_2$. 
Push the $k$ arrows through these arcs with $\Om_2$ and $\Om_5$ moves (see the first step in Figure~\ref{elim_arr}). 
Pushing $k$ arrows through an arc belonging to $C_2$ increases or decreases the sum of signs of crossings between $C_1$ and $C_2$ by $2k$. 
Finally, eliminate all arrows, creating only crossings with both branches in $C_1$. Thus, the linking number between $C_1$ and $C_2$ does not change modulo $k$.
\end{proof}
\end{proposition}

For example, the Hopf link and the trivial $2$-component link are not $2$-equivalent (nor $k$-equivalent for $k\ge 3$).

\begin{figure}[h]
\centering
\includegraphics{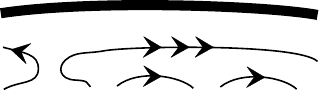}
\caption{$2$ non-exterior and $4$ exterior arrows: $3$ removable and $1$ not}
\label{ext_arr}
\end{figure}

A strand of an arrow diagram that is next to the boundary of the diagram is called {\it exterior}. An arrow is {\it exterior} if it is on an exterior strand.
It is {\it removable} if it is exterior and can be removed with $\Om_\infty$ and an $\Om_4$ (so it is clockwise with respect to the boundary of the diagram).
See Figure~\ref{ext_arr}.

\begin{lemma}\label{lem:1arr}
Suppose that $D$ is an arrow diagram with a single arrow. Suppose this arrow is exterior.
Let $D'$ be $D$ with the arrow removed. Then $D$ and $D'$ are connected with Reidemeister moves.
\begin{proof}
If the arrow is removable apply $\Om_\infty$ and $\Om_4$ to eliminate it. As there are no more arrows it is possible
to drag the arc going next to the boundary of the diagram above the rest of the diagram to get $D'$.

If the arrow is not removable use $\Om_1$ and $\Om_5$ to create a kink and push the arrow into it (see the second step in Figure~\ref{elim_arr}).
Now it is removable. Proceed as in the previous case and, at the end, remove the kink with $\Om_1$ to get $D'$.
\end{proof}
\end{lemma}

To any arrow diagram we associate its {\it shadow} by disregarding the arrows and crossing information.
Two shadows $S$, $S'$ are related by a Reidemeister move $\Om_i$, $1\le i\le 3$, if there are diagrams $D$ and $D'$, 
with respective shadows $S$ and $S'$, where $D'$ is obtained from $D$ with $\Om_i$.

\begin{proposition}
Any knot is $2$-equivalent to the unknot.
\begin{proof}
Let $D$ be an arrow diagram with shadow $S$ and $S'$ a shadow obtained from $S$ with an arbitrary Reidemeister move.
We claim that there exists $D'$ with shadow $S'$, where $D'$ is obtained from $D$ with $2$-moves and Reidemeister moves.

Notice that a $2$-move combined with $\Om_4$ allows to change any arrow into an opposite one.
Let $R$ be a region of $S$ on which we want to perform a Reidemeister move.

If $R$ bounds a loop ($\Om_1$ move), use $2$-moves to reduce the number of arrows on this loop to $0$ or $1$. If there is one arrow
on the loop, push it out of the loop with $\Om_5$, if possible. If not, reverse the arrow first, then push it out. 
As there are no arrows on the loop, apply $\Om_1$ to remove it, getting the required $D'$ with shadow $S'$.

If $R$ bounds a $2$-gon or a triangle ($\Om_2$ or $\Om_3$ move), start by pushing all arrows onto a single exterior strand with $\Om_5$ moves,
reversing the arrows whenever necessary. Then, using $2$-moves, reduce the number of arrows to $0$ or $1$. 
If one arrow remains, use Lemma~\ref{lem:1arr} to eliminate it.
Notice that the shadow $S$ is unchanged.

If $\Om_2$ or $\Om_3$ can be performed on $R$ we get the required $D'$ with shadow $S'$.
Otherwise we want a side of $R$, say $s$, to be exterior. 
Move all arcs separting this side from the boundary of the diagram above the rest of the diagram, so that they are parallel to the boundary
of the diagram, see Figure~\ref{3mv}. We assume that the moved arcs are not involved in any crossing.
Now use Lemma~\ref{lem:1arr} to put an appropriate arrow on $s$, so that a crossing of $R$ can be changed with $\Om_5$ by pushing this arrow.
Now, apply $\Om_2$ or $\Om_3$ on $R$ (it is possible to perform such move, after any single crossing change in $R$). 
Push the arrow next to the boundary with $\Om_5$ moves, reversing it whenever necessary, then eliminate it using Lemma~\ref{lem:1arr}.
Finally, move back the arcs that were changed when making $s$ exterior.
We get the required $D'$ with shadow $S'$.

\begin{figure}[h]
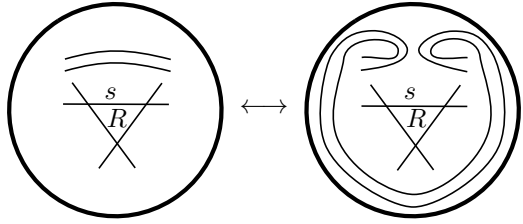

\centering
\oversimage{3mv}{45,55}{$s$}{46,42}{$R$}  \raisebox{4 em}{$\;\longleftrightarrow \;$} \oversimage{3mv2}{46,55}{$s$}{47,42}{$R$}
\caption{Putting $s$ close to the boundary of the diagram}
\label{3mv}
\end{figure}

Obviously, any shadow can be transformed into a circle, if all Reidemeister moves are allowed.
Thus, any diagram $D$ is $2$-equivalent to $D'$ with shadow a circle. Using $2$-moves and $\Om_\infty$, $D'$ can be transformed into a circle with no
arrows on it, i.e. a diagram of the unknot.
\end{proof}
\end{proposition}

\section{Jones polynomial and the $k$-moves}\label{sec:jones}
We use the Kauffman bracket, $<>$, in order to compute the Jones polynomial (see \cite{K}).
Recall that the Kauffman bracket associates to a framed unoriented link an element of $\Z[A,A^{-1}]$. 
Given a diagram of such link with blackboard framing, its Kauffman
bracket is calculated using the following {\it Kauffman relations} (and assuming $<Unknot>=1$):

\[<\raisebox{-7pt}{\includegraphics{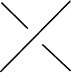}}>=A< \raisebox{-7pt}{\includegraphics{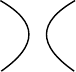}}
>+A^{-1}< \raisebox{-7pt}{\includegraphics{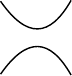}}>\]
\[<L\sqcup\mbox{Unknot}>=(-A^2-A^{-2})<L>\]

For arrow diagrams, it was explained in~\cite{M}, that one can use a framing orthogonal to the fibers of the Hopf fibration, so
that it is preserved under all Reidemeister moves except $\Om_1$. In particular, it is preserved under $\Om_\infty$. 
Such a framing can be deformed into the blackboard framing, which is used from now on.

The Jones polynomial of an oriented link $L$, having a diagram $D$ with writhe $w(D)$, is obtained from the Kauffman bracket via the formula:

\[V_L(t)=(-A)^{-3w(D)}<D>\mbox{, with }t=A^{-4}\]

We will use both $A$ and $t$, with the understanding that $t=A^{-4}$.

Recall, from the previous section, that $T_n$ is the torus knot $T(n,n+1)$, if $n\ge 0$ and it equals $T_{-n-1}$, if $n<0$.
The following formula from~\cite{J} for the Jones polynomials of torus knots is valid, as can be easily checked, not only when $n\in\N$, but
for all $n\in\Z$:

\[V_{T_n}(t)=t^\frac{n(n-1)}{2}\frac{1-t^{n+1}-t^{n+2}+t^{2n+1}}{1-t^2}\label{eq:VT}\tag{VT}\]

It was shown in \cite{M}, that the writhe of $T_n$ is $n(n+1)$, for $n\in\Z$.
Since $(-1)^{3n(n+1)}=1$, it follows that:
\[<T_n>=A^{3n(n+1)}t^\frac{n(n-1)}{2}\frac{1-t^{n+1}-t^{n+2}+t^{2n+1}}{1-t^2}\]

Let $\zeta_k$ be a primitive $k$-th root of unity.

\begin{lemma}\label{lem:Tn_zeta_k}
 Suppose that $k\ge 3$, $n\in\Z$. Then:
\[V_{T_n}(\zeta_k)=-V_{T_{n+k}}(\zeta_k)\mbox{, if $k$ is even}\]
\[V_{T_n}(\zeta_k)=V_{T_{n+k}}(\zeta_k)\mbox{, if $k$ is odd}\]
\begin{proof}
The fraction part of $V_{T_{n+k}}(\zeta_k)$ is
\[\frac{1-\zeta_k^{n+k+1}-\zeta_k^{n+k+2}+\zeta_k^{2(n+k)+1}}{1-\zeta_k^2}
=\frac{1-\zeta_k^{n+1}-\zeta_k^{n+2}+\zeta_k^{2n+1}}{1-\zeta_k^2}\]

The $t$ factor of $V_{T_{n+k}}(\zeta_k)$ is
\[\zeta_k^{\frac{(n+k)(n+k-1)}{2}}=\zeta_k^{\frac{n(n-1)}{2}}\zeta_k^{kn}\zeta_k^{\frac{k(k-1)}{2}}=\zeta_k^{\frac{n(n-1)}{2}}\zeta_k^{\frac{k(k-1)}{2}}\]
If $k$ is even, then $\zeta_k^{\frac{k(k-1)}{2}}=(\zeta_k^\frac{k}{2})^{k-1}=(-1)^{k-1}=-1$.\\
If $k$ is odd, then $\zeta_k^{\frac{k(k-1)}{2}}=1$.
\end{proof}
\end{lemma}

Using $\Om_1$, $\Om_5$ and Kauffman relations one checks easily the following lemma (see~\cite{MD}, Lemmas~3.3 and 3.6):

\begin{lemma}\label{lem:x_inout}
Assume that $2$ (resp. $3$) arrow diagrams are the same except where pictured in the first (resp. second) formula below.
Then their Kauffman brackets are related by the formulas:\\[1.5em]
\raisebox{0.5em}{$<$}\includegraphics{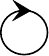}\raisebox{0.5em}{$>\;=A^{-6}\;$}\raisebox{0.5em}{$<$}\includegraphics{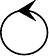}\raisebox{0.5em}{$>$} \\[1.5em]
\raisebox{0.5em}{$<$}\includegraphics{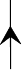}\raisebox{0.5em}{$>\;=-A^2\;<$}\includegraphics{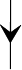}\raisebox{0.5em}{$>\;-A^{-2}\;<$}
\includegraphics{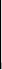}$\;$\includegraphics{x}\raisebox{0.5em}{$>$}
\end{lemma}

We say that a $k$-move is {\it exterior}, if the $k$ arrows involved in it (eliminated or created) are exterior. 
Furthermore, an exterior $k$-move is {\it clockwise} if these $k$ arrows are removable and {\it counterclockwise} otherwise.

\begin{lemma}\label{lem:kbkmove}
Let $D$ be an arrow diagram of a link and let $D'$ be obtained from $D$ by performing a $k$-move on $D$, adding $k$ arrows.
Then the Kauffman brackets of $D$ and $D'$ can be expressed as the following finite sums with the same coefficients $c_i\in\Z[A,A^{-1}]$, $n_i\in\Z$:
\[<D>=\Sum_i c_i <T_{n_i}>\]
\[<D'>=\Sum_i c_i <T_{n_i+k}>\]
Furthermore, if the $k$-move is exterior counterclockwise and $D$ has $n$ arrows, then, in the preceding sums, it is possible to have $n_i\equiv n\mod 2$ for all $i$.
\begin{proof}
By performing on $D$ and $D'$ several $\Om_\infty$ moves if necessary, we can assume that the $k$-move is exterior counterclockwise.
Now, using Kauffman relations, smooth all crossings in $D$ and $D'$ and remove all trivial components, obtaining the same
linear expressions in ovals with arrows for $D$ and $D'$, except for one oval $O$ containing the strand $S$ where the $k$-move is performed.

Let $x$ be an oval with $1$ counterclockwise arrow on it and no ovals nested in it.
Using Lemma~\ref{lem:x_inout} express all ovals, except $O$, with $x$'s: start with expressing the most nested ovals with $x$'s, 
then integrate these $x$'s into other ovals and continue in this way for all ovals except $O$. 
Finally, integrate all $x$'s inside and outside $O$ into $O$. All this can be done in the same way for 
$D$ and $D'$, keeping the strand $S$ freezed. At the end, any element in these linear expressions will be $O$ with some $n_i\in\Z$ arrows on it
for $D$ and with $n_i+k$ arrows for $D'$.

For the last assertion, assume that the $k$-move is counterclockwise and $D$ has $n$ arrows.
One checks easily that transforming $D$ into the linear sum of $T_{n_i}$'s, as above, preserves the number of arrows modulo $2$.
\end{proof}
\end{lemma}

In the next lemma, we assume that the arrow diagram is {\it oriented}, i.e. each of its components has an orientation 
(it is an arrow diagram of an oriented link).

\begin{lemma}\label{lem:wr_kmove}
Let $D$ be an oriented diagram with $n$ arrows.
Let $D'$ be obtained from $D$ by performing an exterior counterclockwise $k$-move, which adds $k$ arrows.
Then:
\[w(D')=w(D)+2kb+k(k+1)\mbox{, for some $b\in\Z$, $b\equiv n\mod 2$}\]
\begin{proof}
Using Reidemeister moves eliminate in the same way all arrows on $D$ and $D'$ (see Figure~\ref{elim_arr}), 
except for the $k$ arrows created on $D'$ by the $k$-move.
This may require some $\Om_1$ moves, but, as they are preformed at the same time on $D$ and $D'$, the difference between the writhes of $D$ and $D'$
is unchanged. Each arrow elimination (with $\Om_\infty$ and $\Om_4$) creates an arc separting the $k$ arrows from the boundary of the diagram.
Thus, we may assume that $D$ has no arrows and $D'$ has $k$ arrows which are separated from the boundary of the diagram by $n$
arcs.

Push each of the $k$ arrows of $D'$ through the $n$ arcs with $\Om_2$ and $\Om_5$ moves. This creates $2kn$ crossings with same signs on a single arc, so
the sum of the signs of these crossings is $2kb$, for some $b\in\Z$, $b\equiv n\mod 2$.

To make the $k$ arrows removable create a negative kink next to each of them with $\Om_1$ and push the arrows into the kinks with $\Om_5$.
The $\Om_1$ moves {\it decrease} the writhe of $D'$ by $k$ and the sum of the new crossings at the end of this step is $k$. Overall it
contributes $2k$ to the difference between the writhes of $D'$ and $D$.

Eliminate one arrow with $\Om_\infty$ and $\Om_4$ and push the other arrows through the new arc creating $2(k-1)$ positive crossings: the
fact that these crossings are positive follows from the assumption that the $k$ arrows are initially on the same strand, next to each other,
see Figure~\ref{wr_plus}.

\begin{figure}[h]
\centering
\includegraphics{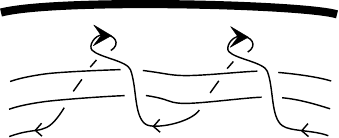}  \raisebox{2 em}{$\;\longrightarrow \;$} \includegraphics{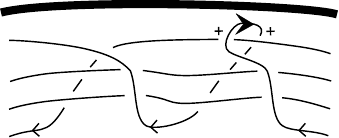}
\caption{Two positive crossings}
\label{wr_plus}
\end{figure}

Repeat this step with the next arrow creating $2(k-2)$ positive crossings. Continue until all arrows are eliminated.

In this way we get:
\[w(D')=w(D)+2kb+2k+2(k-1)+2(k-2)+\ldots+2=w(D)+2kb+k(k+1)\]
\end{proof}
\end{lemma}

In the following theorem we substitute $t=\zeta_k$ in the Jones polynomial of an oriented link $L$. This creates no problems when $L$ is a knot or, more
generally, is a link with an odd number of components. If $L$ has an even number of components, its Jones polynomial has the form $\sqrt{t}$ times
a Laurent polynomial in $t$. In that case we have to choose a square root of $\zeta_k$. 
Therefore, in the following theorem we make the assumption: 
$A$ is a primitive $4k$-th root of unity, $\sqrt{t}=A^{-2}$ and $t=A^{-4}=\zeta_k$.

\begin{theorem}\label{thm:modJones}
Let $D$ be an arrow diagram of an oriented link $L$ and let $D'$ be obtained from $D$ by performing a $k$-move on $D$, $k\ge 3$. Let $L'$
be the link with diagram $D'$. Then:
\[V_L(\zeta_k)=V_{L'}(\zeta_k)\mbox{, if $k$ is odd}\]
\[V_L(\zeta_k)=-V_{L'}(\zeta_k)\mbox{, if $k$ is even}\]
\begin{proof}
By reversing the roles of $D$ and $D'$ if necessary, we may assume that the $k$-move adds $k$ arrows.
As in the proof of Lemma~\ref{lem:kbkmove}, after some $\Om_\infty$ moves, we may also assume that the move is exterior counterclockwise. 

Using the notation of Lemma~\ref{lem:kbkmove} we have:
\[V_D(t)=A^{-3w(D)}<D>=\Sum_i c_i A^{-3w(D)}<T_{n_i}>\]
\[V_{D'}(t)=A^{-3w(D')}<D'>=\Sum_i c_i A^{-3w(D')}<T_{n_i+k}>\]

Let $\epsilon=(-1)^{k+1}$.
It is sufficient to show that, for $t=\zeta_k=A^{-4}$ and any $i$,
\[A^{-3w(D)}<T_{n_i}>=\epsilon A^{-3w(D')}<T_{n_i+k}>\]
or, since $(-1)^{3n_i(n_i+1)}=1=(-1)^{3(n_i+k)(n_i+k+1)}$,
\[A^{-3w(D)}A^{3n_i(n_i+1)}V_{T_{n_i}}(t)=\epsilon A^{-3w(D')}A^{3(n_i+k)(n_i+k+1)}V_{T_{n_i+k}}(t)\]
From Lemma~\ref{lem:Tn_zeta_k},
\[V_{T_{n_i}}(\zeta_k)=\epsilon V_{T_{n_i+k}}(\zeta_k)\]
so one only needs to check that
\[A^{-3w(D)}A^{3n_i(n_i+1)}=A^{-3w(D')}A^{3(n_i+k)(n_i+k+1)}\]
or
\[A^{3(-w(D')+w(D)+(n_i+k)(n_i+k+1)-n_i(n_i+1))}=1\]
Denote by $s$ the exponent of $A$ in the previous equation divided by $3$.
Let $n$ be the number of arrows in $D$. From Lemma~\ref{lem:wr_kmove}:
\[s=-w(D')+w(D)+(n_i+k)(n_i+k+1)-n_i(n_i+1)=\]
\[-2kb-k(k+1)+kn_i+k(k+1)+kn_i=2k(-b+n_i)\]
As  $b\equiv n\mod 2$ and $n_i\equiv n\mod 2$, $-b+n_i$ is even.
Thus $2k(-b+n_i)=4km$, for some $m\in\Z$.
Now 
\[A^{3s}=A^{3(4km)}=t^{-3km}=\zeta_k^{-3km}=1\]
\end{proof}
\end{theorem}

Notice that the case $k=3$ of the theorem is uninteresting as $V_L(\zeta_3)=(-1)^{c-1}$ for any link $L$ with $c$ components, see~\cite{J}.

We now state some corollaries of the previous theorem.

\begin{corollary}\label{cor:main}
Suppose that $L$ and $L'$ are related by a single $k$-move, $k\ge 3$, as in Theorem~\ref{thm:modJones}. Then we have:
\begin{enumerate}[label=(\roman*)]
\item \label{kab} Let $k=ab$, $a\ge 3$. Let $\epsilon=-1$, if $a$ is even and $b$ is odd, and $\epsilon=1$ in all other cases.
Then $V_L(\zeta_a)=\epsilon V_{L'}(\zeta_a)$.
\item \label{k4} If $k=4$, then $L$ and $L'$ either have different Arf invariants; or none of them have an Arf invariant.
\item \label{k8} If $k=8$, then $L$ and $L'$ either have the same Arf invariant; or none of them have an Arf invariant.
\item \label{k6} If $k=6$, then $L$ and $L'$ have the same number of $3$-colorings.
\end{enumerate}
\begin{proof}
\ref{kab}: A $k$-move is obtained by performing $b$ consecutive $a$-moves. Thus, if $a$ is odd, $V_L(\zeta_a)=V_{L'}(\zeta_a)$.
If $a$ is even, then $V_L(\zeta_a)=(-1)^b V_{L'}(\zeta_a)$.

\ref{k4}: With $\zeta_4=i$, $V_L(i)=-V_{L'}(i)$. From~\cite{J}, $V_L(i)=0$ unless $arf(L)$ exists, in which
case $V_L(i)=(-2\sqrt{2})^{c-1}(-1)^{arf(L)}$, where $c$ is the number of components of $L$. Thus, $V_{L'}(i)=V_L(i)=0$
or $arf(L')\neq arf(L)$.

\ref{k8}: it follows from \ref{k4} (an $8$-move is obtained with two $4$-moves).

\ref{k6}: With $\zeta_6=e^\frac{i\pi}{3}$, $V_L(e^\frac{i\pi}{3})=-V_{L'}(e^\frac{i\pi}{3})$.
Thus $3|V^2_L(e^\frac{i\pi}{3})|=3|V^2_{L'}(e^\frac{i\pi}{3})|$.
But these are the numbers of $3$-colorings of $L$ and $L'$, see~\cite{P}.
\end{proof}
\end{corollary}

Notice that the number of $3$-colorings is not preserved under $3$-equivalence: $T_2$ has non trivial 3-colorings but it is $3$-equivalent to the trivial knot.

\begin{corollary}\label{cor:mododd}
Suppose that $L$ and $L'$, with $c$ components, are related by a $k$-move, $k\ge 3$. If $k$ is odd then: 
\[V_L(t)^\equiv V_{L'}(t)\mod (t^k-1)\mbox{, if $c$ is odd}\]
\[\sqrt{t}V_L(t)^\equiv \sqrt{t}V_{L'}(t)\mod (t^k-1)\mbox{, if $c$ is even}\]
\begin{proof}
Let $a>1$ be a divisor of $k$. By Corollary~\ref{cor:main}\ref{kab}, $V_L(\zeta_a)=V_{L'}(\zeta_a)$, as $a$ is odd.
This equality holds also for $a=1$, since $V_L(1)=(-2)^{c-1}$ for any link $L$ with $c$ components, see~\cite{J}.

If $c$ is odd, then $V_L(t)$ and $V_{L'}(t)$ are equal modulo the minimal polynomials (cyclotomic polynomials) of all $\zeta_a$ where $a$ divides $k$.
It follows that they are equal modulo the product of these polynomials, which is $t^k-1$.

If $c$ is even, the Laurent polynomials $\sqrt{t}V_L(t)$ and $\sqrt{t}V_{L'}(t)$ are equal in $\zeta_a$ for all divisors $a$ of $k$
(we choose a square root of each $\zeta_a$, see the discussion before Theorem~\ref{thm:modJones}). So, again, they are equal modulo $t^k-1$.
\end{proof}
\end{corollary}

\begin{figure}
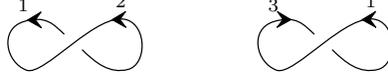
\centering
\oversimage{revK_1}{8,44}{\scriptsize $1$}{80,47}{\scriptsize $2$}
\qquad\qquad \oversimage{revert_1}{8,44}{\scriptsize $3$}{80,47}{\scriptsize $1$}
\caption{Knots $9_{42}$ and $10_{124}$}
\label{5-move}
\end{figure}

As an example, consider the knots with arrow diagrams shown in Figure~\ref{5-move}. These knots were considered in~\cite{M}, the first is 
$9_{42}$ and the second $10_{124}$ (or the torus knot $T(3,5))$. 
Rotate the first diagram by $\pi$ and apply a $5$-move by adding $5$ arrows opposite to the $2$ arrows to
get the second diagram. Thus, these knots are related by a single $5$-move. Their Jones polynomials are: 

\[V_{9_{42}}(t)=t^{-3}-t^{-2}+t^{-1}-1+t-t^2+t^3\equiv -1+t+t^4\mod (t^5-1)\]
\[V_{10_{124}}(t)=t^4+t^6-t^{10}\equiv -1+t+t^4\mod (t^5-1)\]

It follows, that these two knots are not $5$-equivalent to the unknot.

One can separate many $5$-equivalence classes for knots, say, up to $10$ crossings, by checking their Jones polynomials modulo $t^5-1$.
For example, the knot $7_7$ is not $5$-equivalent to any other knot with at most $10$ crossings (it may be equivalent to at most $3$ knots
with $11$ crossings and $4$ knots with $12$-crossings). Of course, equality modulo $t^5-1$ of the Jones polynomials of two knots does not guarantee that 
these knots are $5$-equivalent.

When $k$ is even, a version of Corollary~\ref{cor:mododd} is more complicated because of the change of sign of the Jones polynomial in $\zeta_k$ 
and the fact that Theorem~\ref{thm:modJones} does not hold for $k=2$. We have:

\begin{corollary}\label{cor:modeven}
Suppose that $L$ and $L'$, with $c$ components, are related by a $k$-move, $k\ge 3$. If $k$ is even then:
\[V_L(t)^\equiv t^\frac{k}{2}V_{L'}(t)\mod \left(\frac{t^k-1}{t+1}\right)\mbox{, if $c$ is odd}\]
\[\sqrt{t}V_L(t)^\equiv t^\frac{k}{2}\sqrt{t}V_{L'}(t)\mod \left(\frac{t^k-1}{t+1}\right)\mbox{, if $c$ is even}\]
\begin{proof}
Let $k=ab$, $a\neq 2$. If $b$ is even and $a\ge 3$, then $V_L(\zeta_a)=V_{L'}(\zeta_a)$ from Corollary~\ref{cor:main}\ref{kab}. Such equality holds
also for $a=1$, as in the proof of Corollary~\ref{cor:mododd}.
Also, if $b$ is even, $\zeta_a^\frac{k}{2}=\zeta_a^\frac{ab}{2}=1$. 
If $b$ is odd, so that $a$ is even, then $V_L(\zeta_a)=-V_{L'}(\zeta_a)$ from Corollary~\ref{cor:main}\ref{kab}. 
Also, if $b$ is odd, $\zeta_a^\frac{k}{2}=\zeta_a^\frac{ab}{2}=(-1)^b=-1$.

If $c$ is odd, then, for any divisor $a$ of $k$, $a\neq 2$, $V_L(t)$ and $t^\frac{k}{2}V_{L'}(t)$ are equal modulo the minimal polynomial of $\zeta_a$.
This means that these two polynomials are equal modulo the product of these cyclotomic polynomials, which is $\frac{t^k-1}{t+1}$, $t+1$ being the minimal polynomial
of $-1$.

For $c$ even, the argument is the same as in the proof of Corollary~\ref{cor:mododd}.
\end{proof}
\end{corollary}

We now study the influence of the $k$-moves on the Jones polynomial in $-1$ (i.e. the determinant up to sign). 
Checking $V_{T_n}(-1)$ one notices that if $k=2p$ then:
\[V_{T_{n+k}}(-1)\equiv (-1)^p V_{T_n}(-1)\mod k\]
This can be used to prove that the same is true for any $L$ and $L'$ related by a $k$-move i.e.
\[V_{L'}(-1)\equiv (-1)^p V_L(-1)\mod k\]
One possible proof follows closely that of Theorem~\ref{thm:modJones}: an equality up to sign has to be replaced by a congruence modulo $2p$.

However, it turns out that the congruence above is a consequence of Corollary~\ref{cor:modeven}. In particular, if this corollary fails to distinguish
different classes modulo $k$, then the congurence will not distinguish them either.

\begin{proposition}\label{pro:detmodp}
Suppose that $k=2p$. Let $f(n)=(n+1)\mod 2$ ($f(n)$ is $0$ or $1$).
Let $L$ and $L'$ be two links, with $c$ components, for which the conclusion of Corollary~\ref{cor:modeven} is true, i.e.:
\[\sqrt{t}^{f(c)}V_L(t)\equiv t^\frac{k}{2}\sqrt{t}^{f(c)}V_{L'}(t)\mod \left(\frac{t^k-1}{t+1}\right)\]
Then
\[i^{f(c)}V_{L'}(-1)\equiv (-1)^p i^{f(c)}V_L(-1)\mod k\]
In particular, if $L$ and $L'$ are knots, $p$ is an odd prime and $L$ is $p$-colorable, then so is $L'$.
\begin{proof}
It suffices to evaluate the congurence in $t=-1$. Notice that
\[\frac{t^k-1}{t+1}=t^{k-1}-t^{k-2}+\ldots+t-1\]
For $t=-1$ it equals $-k$ and $t^\frac{k}{2}$ equals $(-1)^p$.
If $c$ is even, we choose $i$ as the square root of $-1$.

For the last assertion, the $p$-colorability is equivalent to the determinant being divisible by $p$ (see for instance~\cite{Lv}). 
If $V_L(-1)$ is divisible by $p$, then so is $V_{L'}(-1)=-V_L(-1)+2pa$, for some $a\in\Z$.
\end{proof}
\end{proposition}

\section{Applications}\label{sec:appl}

\begin{proposition}\label{pro:tn}
Let $k\ge 5$ and $s=\lfloor\frac{k-1}{2}\rfloor$.
Then for torus knots $T_n$, $n\in\N$, there are $s$ distinct classes modulo $k$ with representatives $T_1,T_2,\ldots,T_s$.
\begin{proof}
Using $k$-moves it is obvious that any $T_n$ is in one of these classes.
We use Theorem~\ref{thm:modJones} to show that they are distinct. In fact, we show that the modules of the Jones polynomials of $T_n$ in $\zeta_k$
are distinct for $n=1,\ldots,s$.

Using equation (\ref{eq:VT}):

$|V_{T_n}(\zeta_k)|=\frac{|1-\zeta_k^{n+1}-\zeta_k^{n+2}+\zeta_k^{2n+1}|}{|1-\zeta_k^2|}$

We can disregard the denominator which does not depend on $n$. Denote the numerator by $N(n)$.
Factoring out $\zeta_k^{(2n+1)/2}$ we get:
\[N(n)=\left|2cos\frac{2n\pi+\pi}{k}-2e^\frac{2i\pi}{k}cos\frac{\pi}{k}\right|\]
So that:
\[N^2(n)=\left(2cos\frac{2n\pi+\pi}{k}-2cos\frac{2\pi}{k}cos\frac{\pi}{k}\right)^2+\left(2sin\frac{2\pi}{k}cos\frac{\pi}{k}\right)^2\]
The derivative of $N^2(n)$ is:
\[\frac{16\pi}{k}\left(cos\frac{2\pi}{k}cos\frac{\pi}{k}-cos\frac{2\pi n+\pi}{k}\right)sin\frac{2n\pi+\pi}{k}\]
We check the sign of this derivative for $n\in [1,s]$.
The second term with the sine is positive for $n\in [1,s)$ (it is zero for $n=\frac{k-1}{2}\ge s$).
For $n=1$, expanding the second cosine one checks that the first term equals $sin\frac{2\pi}{k}sin\frac{\pi}{k}$ which is positive.
Thus the derivative is positive in $n=1$.
Increasing $n$ from $1$, $cos\frac{2\pi n+\pi}{k}$ decreases until it reaches $-1$ in $n=\frac{k-1}{2}\ge s$.
Thus the first term is positive for $n\in [1,s]$.
It follows that $N^2(n)$ is increasing in this interval. In particular it takes different values in $1,2,\ldots,s$ and we are done.
\end{proof}
\end{proposition}

If $k\le 4$ one checks easily that all torus knots $T_n$ are $k$-equivalent to the trivial knot. 
A natural question arises: is every knot $k$-equivalent to the trivial knot when $k=3$ or $4$?

We now turn to applications involving the Hopf crossing number of a link $L$, denoted $h(L)$ (see the Introduction).

\begin{lemma}\label{lem:hopf_fin}
Let $\mathcal L_{c,s}$, $c\ge 0$, $s\ge 1$, be the set of all links $L$ with at most $s$ components satisfying $h(L)\le c$. 
For $k\ge 1$ there are finitely many links in $\mathcal L_c$ up to $k$-equivalence.

In particular, for $k\ge 3$, the set $\{V_L(\zeta_k),L\in\mathcal L_{c,s}\}$ is finite. 
\begin{proof}
Let $D$ be an arrow diagram of a link in $\mathcal L_{c,s}$. Using $k$-moves, one may reduce the number of arrows on each arc to
a number between $0$ and $k-1$. As there are finitely many classical diagrams with at most $c$ crossings and at most $s$ components, it is clear that
there are also finitely many links in $\mathcal L_{c,s}$ up to $k$-equivalence.

The second part follows directly from Theorem~\ref{thm:modJones}.
\end{proof}
\end{lemma}

This lemma can be used to construct knots with arbitrarily high Hopf crossing number. Also, by computing the finite set of values of $V_L(\zeta_k)$,
for a given $k\ge 5$ and $h(L)\le c$, successively for $c=0,1,2,\ldots$, one may find lower bounds of $h$ for concrete knots 
(for instance from Rolfsen's table~\cite{R}).

\begin{figure}[h]
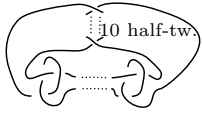

\centering
\overimage{hopf_pa}{48,37}{\scriptsize $10$ half-tw.}
\caption{Family of knots $K'_n$}
\label{hopf_pa}
\end{figure}

Let us consider some examples of knot families with unbounded Hopf crossing number. Denote by $K_n$ the connected sum of $n$ copies of the right
handed trefoil knot (i.e. $T_2$). Let $K'_n$ be the family of knots obtained from $K_n$ be adding $10$ half-twists, as shown in Figure~\ref{hopf_pa}.
Notice that $K'_n$ are alternating and prime, because their diagrams, shown in that figure, are prime (see \cite{L}).

\begin{proposition}
For $K_n$ and $K'_n$ one has: 
\[\lim_{n\to\infty} h(K_n)=\infty\mbox{ and }\lim_{n\to\infty} h(K'_n)=\infty\]
Also, for $k\ge 5$, there are inifinitely many classes of knots modulo $k$.
\begin{proof}
Let $k\ge 5$. It follows from Proposition~\ref{pro:tn}, that $|V_{T_2}(\zeta_k)|>1$.
Now $V_{K_n}(\zeta_k)=(V_{T_2}(\zeta_k))^n$, so the sequence $|V_{K_n}(\zeta_k)|$, $n\in\N$, is increasing.
It follows from Theorem~\ref{thm:modJones}, that $K_n$ and $K_m$ are not $k$-equivalent if $n\neq m$.
From Lemma~\ref{lem:hopf_fin}, we get:
\[\lim_{n\to\infty} h(K_n)=\infty\]
For the family $K'_n$, one notices that it is obtained from $K_n$ by a $\bar{t}_{10}$ move, see~\cite{Ptk}.
It is shown there, that $V_{\bar{t}_{2k}(L)}(t)=V_L(t)$, for $t^{2k}=1$. Taking $k=5$, we get
$|V_{K'_n}(\zeta_{10})|=|V_{K_n}(\zeta_{10})|$. One concludes by applying Lemma~\ref{lem:hopf_fin} again.
\end{proof}
\end{proposition}

Here is another application of Theorem~\ref{thm:modJones}:

\begin{proposition}
Let $\mathcal K_3$ be the family of prime knots having braid index equal to $3$.
The Hopf crossing number is unbounded on $\mathcal K_3$.
\begin{proof}
Let $\mathcal L_3$ be the family of links having braid index $3$.
It follows from~\cite{J}, that $S:=\{|V_L(\zeta_{14})|,L\in\mathcal L_3\}$ is dense in the interval $[0,(2\cos\frac{\pi}{14})^2]$.
On the other hand, any 3-braid, $B$ can be easily transformed into another 3-braid $B'$ with $t_7$ moves (adding $7$ half-twists),
so that the closure of $B'$ is a knot. It follows from~\cite{Ptk}, that the module of the Jones polynomial in $\zeta_{14}$ is
unchanged by $t_7$ moves. Thus $S$ is unchanged if we limit ourselves to knots.

The braid index satisfies $b(K\sharp K')=b(K)+b(K')-1$ (see~\cite{BM}). Thus, the composite knots with braid index $3$
are of the form $T(2,a)\sharp T(2,b)$, for some torus knots $T(2,a)$ and $T(2,b)$. One checks easily that the Jones polynomial
of $T(2,a)$ takes finally many values in $\zeta_{14}$ for $a\in\N$, so the same is true for $T(2,a)\sharp T(2,b)$, $a,b\in\N$.
Thus $S':=\{|V_L(\zeta_{14})|,L\in\mathcal K_3\}$ equals $S$ minus a finite set ($S'$ may be equal to $S$).

As $S'$ is infinite, we conclude by applying Theorem~\ref{thm:modJones}.
\end{proof}
\end{proposition}

It is quite likely that the Hopf crossing number is unbounded for torus knots $T(2,a)$. However, it cannot be proved using Theorem~\ref{thm:modJones}.

We conclude with an interpretation of $k$-moves involving links in lens spaces $L_{p,1}$. For $p\in\Z$, consider the Reidemeister move $\Om^p_\infty$,
which looks like $\Om_\infty$ in Figure~\ref{reid_moves}, except that there are $p$ arrows on the right diagram, instead of $1$ (so that $\Om^1_\infty$ 
is the same as $\Om_\infty$). Taking the same arrow diagrams as for $S^3$ and replacing $\Om_\infty$ with $\Om^p_\infty$, we obtain diagrams and 
Reidemeister moves for links in $L_{p,1}$, see~\cite{M2}. In fact, the convention used in that paper as to the direction of $p>0$ arrows is opposite to
the one used in the present paper, but both choices give Reidemeister moves for $L_{p,1}$, since it is homeomorphic to $L_{-p,1}$. 
In particular, we can define $k$-moves for links in $L_{p,1}$.

Consider $\Om_\infty$, $\Om^{k+1}_\infty$ and $k$-moves. One checks without much difficutly that with any two of these moves one can genarate the third
(for generating a $k$-move on a strand $s$ with the other two, use some $\Om_\infty$ so that $s$ is next to the boundary, apply $\Om_\infty$
and $\Om^{k+1}_\infty$ on $s$ to obtain $k$ arrows, then undo the first step with some $\Om_\infty$).
Thus, the equivalence classes of links in $S^3$ modulo $k$-moves are obtained by formally identifying diagrams of links in $S^3=L_{1,1}$ and $L_{k+1,1}$.
For example, identifying in that way links in $S^3$ and $L_{6,1}$, we obtain equivalence classes modulo $5$-moves.

It is not necessary to start with $L_{1,1}=S^3$. Let $a\in\Z$. Using the same argument as above, one shows that any two moves generate the third in the triple
$\Om^{a}_\infty$, $\Om^{a+k}_\infty$ and $k$-moves. In fact, $\Om^a_\infty$ and $k$-moves generate $\Om^{a+kb}_\infty$ moves, for any $b\in\Z$.
Thus, allowing $k$-moves in $L_{a,1}$ corresponds to identifying links via their diagrams in $L_{a+kb,1}$, $b\in\Z$. Here, one has to be a bit careful: even
though $L_{p,1}\cong L_{-p,1}$, the sign is important in identifying diagrams. For example, an oval with $2$ counterclockwise arrows (or $T_2$) is
a trefoil in $L_{1,1}$, whereas it is a trivial knot in $L_{-1,1}$ (it is the same as an oval with $2$ clockwise arrows in $L_{1,1}$).
To illustrate these identifications, let $k=3$: $L_{0,1}$ is identified with $L_{3b,1}$, $L_{1,1}$ with $L_{1+3b,1}$ and $L_{2,1}$ with $L_{2+3b,1}$, $b\in\Z$.
Although $L_{1,1}$ is identified with $L_{-2,1}$ it is not identified with $L_{2,1}$ (that would generate $1$-moves).

An observation following from Theorem~\ref{thm:modJones} is that the Jones polynomial at $\zeta_k$ is, by such identifications, naturally defined 
(up to sign if $k$ is even) for classes of links in $L_{1+ka,1}$ modulo $k$-moves, $a\in\Z$. For example, it is defined in $\zeta_5$ for links
in $L_{6,1}$ modulo $5$-moves. Finally, it is quite clear that one can define the Hopf crossing number for links in $L_{p,1}$ and use 
Theorem~\ref{thm:modJones} to study it (showing, for example, that there are knots in $L_{6,1}$ with aribtrarily high Hopf crossing number).

\end{document}